\documentclass[12pt,oneside]{article}
\usepackage{amsmath,amssymb,amsfonts,amsthm}
\usepackage{amscd}
\usepackage{graphicx}
\usepackage[all]{xy}
\usepackage{multirow}
\usepackage{rotating}
\usepackage{hyperref}
\usepackage{hyperref}
\hypersetup{
    colorlinks=true,
    linkcolor=blue,
    filecolor=magenta,
    urlcolor=cyan,
    citecolor=blue}
\usepackage{lscape}
\newtheorem{Thm}{Theorem}[section]
\newtheorem{Lem}[Thm]{Lemma}
\newtheorem{Pro}[Thm]{Proposition}

\newtheorem{Def}[Thm]{Definition}
\newtheorem{Rem}[Thm]{Remark}

\theoremstyle{definition}

\theoremstyle{remark}

\renewcommand{\phi}{\varphi}

\numberwithin{equation}{section}

\begin{document}
\title{On the generalized $m$-Kropina metrics}
\author{Ebtsam  H. Taha}
\date{}
\maketitle    
\begin{center}
Department of Mathematics, Faculty of Science, Cairo University, Giza, Egypt,\\
E-mail: ebtsam.taha@sci.cu.edu.eg, ebtsam.h.taha@hotmail.com
\end{center} 
 \begin{center}
 \textit{To the memory of Professor Nabil L. Youssef}
 \end{center}

\maketitle 

\vspace{0.5cm}
\textbf{Abstract.}
Generalized $m$-Kropina metrics appear naturally as a spacetime geometry compatible with  Lorentz symmetry breaking, leading to  useful applications in modified gravity and cosmology. We prove that a generalized $m$-Kropina metric $F$ is an almost rational Finsler metric. Thereby, we study the rationality of its Finslerian geometric objects in the directional variable $y$.   For example, its geodesic spray coefficients are rational in $y$.  Consequently, we prove that if $F$ is an Einstein metric with $m \notin \mathbb{Z}$, then it is Ricci-flat. Moreover, for $m \in 2 \mathbb{Z}$,  the arithmetic nature of $m$ imposes strong rigidity constraints:  if $F$ has isotropic mean Berwald curvature, or has relatively isotropic Landsberg curvature, or has almost vanishing $\mathbf{H}$-curvature, then $F$ is weakly Berwaldian, or $F$ is Landsbergian, or $\mathbf{H}=0$, respectively. Furthermore, we show that if $F$ has almost isotropic flag curvature ($m \in 2 \mathbb{Z}$ and $n \geq 3$), then the flag curvature is constant. We, hence, deduce under what conditions a generalized $m$-Kropina metric $F$ becomes an exact solution to \lq \lq Pfeifer and Wohlfarth’s vacuum field equation".  Finally, we provide several four-dimensional examples arising in modified gravity and cosmology.

\bigskip 

\maketitle 


\maketitle 

\section{Introduction} 
  Kropina metrics belong to the family of $(\alpha ,\beta)$- metrics which are considered as a perturbation of a Riemannian metric $\alpha$ by $1$-form $\beta$.   In fact, $(\alpha ,\beta)$- metrics originally defined by Matsumoto,  in \cite{Matsumoto},  which  include Randers,  Matsumoto, Randers-Kropina and square metrics cf. \cite{Erasmo 23, Shenbook16,  Amr-Taha-Kropina, arFinslerTaha, gen. Krop. change}. For instance, Randers metrics represent solutions to  Zermelo’s  navigation problem in the case of weak wind (see, e.g. \cite{Erasmo} for a wide perspective of Randers or the more general notion of a wind Finsler metric).  Pseudo-Finsler metrics represent a class of anisotropic spacetime models (as its Finsler metric tensor depends not only on position but also on direction) which can help to describe spin, torsion or other non-Riemannian effects \cite{Erasmo 18, Erasmo 23, Hohmann_2019, Pfeifer:2011xi}.  Finsler spacetimes (which are  pseudo-Finsler metrics with certain properties) possess essential characteristics similar to Lorentzian spacetimes, making them suitable for physical applications such as a well-defined causal structure, a clearly defined concept of proper time (which is arc length) as well as Finslerian geodesic equations that can be understood as the trajectories of freely falling particles \cite{Sjors thesis}. 

\bigskip 

\par An m-Kropina metric (which is a generalization of a Kropina metric)  arises naturally as a spacetime geometry consistent with the violation of Lorentz symmetry \cite{Sjors thesis}.  A generalized m-Kropina metric (which is a further generalization of m-Kropina metric) is a possible candidate to be a model of local anisotropy in spacetime, that can be linked to fundamental particles or fields.  It can describe particles which move in a preferred direction (i.e., breaking isotropy). Moreover,  it can allow anisotropic corrections to general theory of relativity which  could have a role in a better understanding of   quantum gravity,  dark matter,  high-energy physics or cosmology  \cite{Fuster:2018djw,  Sjors thesis, Sjors m-Kropina}. The aim of this work is to determine under what conditions they provide be exact solutions to a Finslerian version of Einstein's field equations derived by Pfeifer and Wohlfarth \cite{Hohmann_2019, Pfeifer:2011xi}.
\bigskip

This paper is organized as follows.  In Section 2, we recall some necessary materials to make this work self-contained.  In \S 3,
we give a possible meaning to the parameters  of the generalized $m$-Kropina metric $F$ expression in order to bridge the  mathematical framework with the possible physical applications. We give the conditions under which $F$ is either positive definite, non-degenerate or has Lorentz signature. We prove that a generalized $m$-Kropina metric $F$ is an AR-Finsler metric which facilitates the investigation of the rationality of the associated Finslerian  objects  in $y$. Then, we prove that a weakly Einstein metric $F$, $m \in 2\mathbb{Z}$, is  Einstein. Further, an Einstein metric $F$ with $m \notin \mathbb{Z}$, is Ricci-flat. Moreover, for $m \in 2 \mathbb{Z}$, if $F$ has isotropic mean Berwald curvature, or has relatively isotropic Landsberg curvature, or has almost vanishing $\mathbf{H}$-curvature, then $F$ is weakly Berwaldian, or $F$ is Landsbergian, or $\mathbf{H}=0$, respectively. Also, we show that if $F$ has almost isotropic flag curvature ($m \in 2 \mathbb{Z}$ and $n \geq 3$), then the flag curvature is constant.  We end this section by a table which is created to summarize the rationality of the Finslerian geometric objects related to $F$.  In \S 4 we make use of our results to provide some conditions on the family of generalized $m$-Kropina metrics ensuring that they satisfy  \lq \lq Pfeifer and Wohlfarth’s vacuum field equation".  Also,  we  give some well-known examples which arise from the family of generalized $m$-Kropina metrics,  from the Very Special Relativity (VSR),  disformal gravity to Bumblebee and Einstein-Aether gravity models.  Finally, we summarize our results along with future work in the conclusion.
\section{Basics of pseudo-Finsler spaces}
We start with the following notation: $M$ is an $n$-dimensional, $n>1$, smooth manifold with points $x=:( x^1 , x^2, ..., x^n )$.  Also, $C^\infty (M)$ denotes the space of smooth functions on $M$.  $TM$ denotes the tangent bundle with fibres $T_x M$ and points $(x,y)$ such that $x \in M$ and $( y^1 , y^2, ..., y^n )=:y \in T_{x}M$. An open subset $\mathcal{A}$ of $TM  \setminus \{0\}$ is conic if it satisfies for every $(x,y) \in \mathcal{A}$ and $\mu  \in \mathbb{R}^+ ,\, (x,\mu\,y) \in  \mathcal{A} $. A conic subbundle $\mathcal{A}$ of $T M  \setminus\{0\}$, is a conic open subset of $T M$ such that its projection $\pi (\mathcal{A})=M$.   $\mathcal{A}$ has  fibres denoted by $\mathcal{A}_{x}$. 
Let us denote the differentiation with respect to $x$ by $\partial_i$ and the differentiation with respect to $y$ by $\dot{\partial}_i$.

 \begin{Def}\label{pseudo-Finsler def}
  A conic pseudo-Finsler structure \footnote{Several stricter definitions of Finsler spacetimes exist in the literature, starting from Beem's original definition \cite{Beem}. Their specific technical requirements differ based on the application's scope; for example, look to \cite{Erasmo 20, Voicu 22, Finsler definitions, ref5} and the references therein.}  $ F$ on $ M $ is a smooth  function   $\;F:\mathcal{A} \longrightarrow \mathbb{R}\,$\footnote{For the generalized m-Kropina metrics considered below, the conic domain $\mathcal{A}$  will be specified explicitly.}   such that:
  \begin{description}
    \item[(i)]  $ F$ is positively homogeneous of one degree  in $y$:  $F(x,\lambda y)=\lambda \,F(x,y),\,\,\,\forall  \lambda \in \mathbb{R^+}.$
    \item[(ii)] The Finsler metric tensor, whose components $ g_{ij} (x,y):=\frac{1}{2}\dot{\partial}_{i}\dot{\partial}_{j}F^2(x,y),$  is non-degenerate at each point of  $ \mathcal{A}_{x}, \,\,\forall x \in M$. 
  \end{description} \end{Def} 
 Since some results we present here hold for Finsler spaces with positive definite metric, for Finsler spacetimes with Lorentzian signature and for any other choice of signature of the Finsler metric one may like to put in place, we do not specify the Finsler metric's signature. The Finsler metric just needs to be invertible on a sufficiently large subset of the tangent bundle.\\ 

Associated to a Finsler metric $F$, many non-Riemannian geometric objects are defined such as the components of Cartan tensor $C_{ijk}=\frac{1}{2} \dot{\partial}_{k}\, g_{ij}$ and $h_{ij} :=F \dot{\partial}_i \dot{\partial}_j F $ are the angular metric components.  
In addition,  the Finslerian geodesic spray which has coefficients $G^i= \frac{1}{4} \, g^{ir}\left( y^k \, \dot{\partial}_r \partial_k F^2 - \partial_r   F^2 \right),$  determined by the goedeic equations \cite{ChernShen2005}:
\[\frac{d^2 x^i (t)}{d t^2} + G^i \left(x (t) ,\frac{d x (t)}{d t} \right) =0,\quad i=1,...,n. \] 

Consequently, the Cartan nonlinear (Barthel) connection coefficients are given by $N^i_j := \dot{\partial}_{j} G^i $. Once the Barthel connection is defined, the double tangent bundle splits into the vertical bundle (with basis $\{\dot{\partial}_j \}^{n}_{j=1}$) and horizontal bundle (which spanned by $\{\delta _i := \partial _i -N^{j}_{i}\, \dot{\partial}_j\}^{n}_{i=1}$).  The curvature of Barthel connection has components $R^{i}_{jk}$ defined as follows \begin{equation}\label{curvature of Barthel connection}
R^{i}_{jk} := \delta_{k} N^{i}_{j} - \delta_{j} N^{i}_{k} = \partial_{k} N^{i}_{j} - \partial_{j} N^{i}_{k} + N^{r}_{k} G^{i}_{rj} -N^{r}_{j} G^{i}_{rk}.
\end{equation}
Thus, the Finsler Ricci scalar (or, simply Ricci scalar) is defined by $Ric= R^{k}_{ij}\, y^j$. The Ricci scalar also given by $Ric=R^i_i $, where the Riemannian curvature components $R^i_k$ are given by \cite{ChernShen2005}
\begin{equation}\label{h Riem curv}
R^i_k :=2\,\partial_k G^i - y^{j}\,\partial_j\,\dot{\partial}_k G^i +2 G^j\,\dot{\partial}_k \,\dot{\partial}_j G^i - \dot{\partial}_j G^i\, \dot{\partial}_k G^j .
\end{equation} 
Moreover, the components of the Finsler Ricci tensor is given by $R_{ij}=\frac{1}{2} \dot{\partial}_{i} \dot{\partial}_{j}  Ric$.
It is known that, the vanishing of $Ric$ is equivalent to the vanishing of $R_{ij}$ which follows from 
\begin{equation}\label{Ricci scalar and tensor} 
R_{ij}=\frac{1}{2} \dot{\partial}_{i} \dot{\partial}_{j}  Ric \quad \text{ and } \quad Ric= R_{ij}y^i y^j.
\end{equation}
\par The Berwald curvature components are defined by $G^{i}_{jkl}=\dot{\partial}_{j}\dot{\partial}_{k}\dot{\partial}_{l} G^i$. Moreover,  the Landsberg curvature components are defined by $L_{ijk}=-\frac{1}{4}\,y^{l}g_{ls}\,G^s_{ijk}$ and mean Landsberg curvature is given by $J :=J_k\, dx^k = g^{ij}L_{ijk} \, dx^k. $ It is known that a Finsler metric $F$ is called Berwaldian if $G^{i}_{jkl}=0$. Also, $F$ is called weakly Berwaldian if the mean Berwald curvature $E_{ij} := \frac{1}{2} G^{k}_{kij} $ vanishes. On the other hand, $F$ is called  Landsbergian  if $L_{jkl}=0$ and is called weakly Landsbergian if $J =0$ \cite{Sjors thesis}.
\begin{Def}\emph{\cite{Shenbook16, ChernShen2005} }\label{Finsler metric types}
		A Finsler metric $F$ on a manifold $M$ is said to be of
\begin{description}
			\item[(a)]  isotropic mean Berwald type if there exists a smooth  function $A$ on $M$ such that $$E_{ij}= \frac{(n+1)}{2 F}\, A \,h_{ij}.$$
			\item[(b)]  isotropic mean Landsberg type if $J$ can be expressed as $J = A\, F I$ for some smooth function $A$ on $M$.
			\item[(c)] relatively isotropic Landsberg type if there exists a smooth  function $A$ on $M$ such that $L_{ijk}= A\, F\,  C_{ijk}$.
			\item[(d)]  isotropic $S_{\mu}$-curvature  type if its $S_{\mu}$-curvature associated with an arbitrary volume form $d V_{\mu}= \sigma_{\mu} (x) dx$ can be written as $S_{\mu}= (n+1)  F\, A $ for some smooth  function $A $ on $M$, where $S_{\mu}(x,y):=N^i_i (x, y)-y^i\, \partial_{i}\log\left( \sigma_{\mu}(x)\right)$ and $ \sigma_{\mu}$ is the volume density function of an arbitrary volume measure $d V_\mu$ on M.
			\item[(e)] Weak Einstein type if the Finsler Ricci scalar $Ric$ can be written as $$Ric = (n-1) \{ \frac{3 \theta}{F} + K \}F^2 ,$$  for some smooth  function $K$ on $M$ and a $1$-form $\theta$ on $M$.  In particular,  $F$ is an Einstein metric if $\theta =0$, or equivalently, $Ric = (n-1) K F^2$.
			\item[(f)] Ricci-flat type if its Finsler Ricci-tensor vanishes identically.
			\item[(g)] almost isotropic flag curvature $\mathbf{K}$ if $\mathbf{K}$  can be written as $\mathbf{K}= \frac{y^i \partial_{i} A}{F} + \sigma $, for some smooth functions  $A \text{ and } \sigma$ on $M$.
			\item[(h)] constant flag curvature  if  $\mathbf{K}$ is constant.
			\item[(i)] almost vanishing $\textbf{H}$	if there exists a smooth function $\rho$ on $M$ such that $$H_{ij} := y^l E_{ij||l}=y^l ( \delta _{l} E_{ij} - G^{k}_{il} E_{kj}- G^{k}_{jl} E_{ki})$$ satisfies $F H_{ij}= \frac{(n+1)}{2 } A \, h_{ij}$, where $(||)$ is the horizontal  covariant derivative with respect to Berwald connection, for some smooth functions  $A$ on $M$. 
			\end{description}	\end{Def}	
			 \section{Analysis of Generalized $m$-Kropina metrics }	
			 An interesting family of Finsler metrics is the one of $(\alpha ,\beta)$-metrics which is introduced by Matsumoto in \cite{Matsumoto}. Let $\alpha = \sqrt{|\alpha_{ij}(x) y^iy^j|}$ with $\alpha_{ij}(x)$ be  a pseudo-Riemannian metric on $M$  and $\beta = b_i(x) y^i$ a $1$-form on $M$. Assume that $\alpha\neq 0$ on $\mathcal A$. A Finsler metric $F $ is said to be an indefinite $(\alpha,\beta)$-metric, with domain $\mathcal A$ not including $\alpha=0$, if it can be written as $F = \alpha\phi(s) $, where $\phi$ is smooth on its domain and  $s := \beta/\alpha$ viewed as real variable of $\phi$ \cite{Sjors thesis}. Since $\alpha$ is not necessarily a positive definite Riemannian metric, we have  $\alpha ^2 = sgn (\alpha ^2)\, \alpha_{ij}(x) y^i y^j $, where $sgn(\alpha ^2)$ is the sign of $\alpha ^2$. Thus, \[\alpha_{i}:=\dot{\partial}_{i}\alpha= \frac{sgn (\alpha ^2) \alpha_{ij}\, y^{j}}{\alpha} =sgn (\alpha ^2)\, \frac{ y_{i}}{\alpha}, \quad y_{i}:= y^j\, \alpha_{ij}, \quad  \alpha_{i} \, \alpha^{i} =sgn (\alpha ^2).\footnote{Unless stated, we raise and lower the indicies using the pesudo Riemannian metric  $\alpha _{ij}$}\]
	Consequently, the metric tensor components for a pseudo-Finsler metric of $(\alpha , \beta)$-type (or, an indefinite $(\alpha,\beta)$-metric) are given by \cite{ref5}
			\begin{eqnarray} \label{gen metric}
			g_{ij}&=& \left(\phi^2 -s \phi\, \phi' \right)\,sgn (\alpha ^2)\,  \alpha_{ij} + \left( \phi \, \phi'' + (\phi')^{2}\right)\, b_{i}\, b_{j} \nonumber\\ 
			& &+ \frac{1}{\alpha}\left(\phi\, \phi' -s\left[\phi\,\phi'' + (\phi')^{2} \right] \right) \, \{sgn (\alpha ^2)\,( b_{i} y_{j} +  b_{j} y_{i}) -  \frac{s}{\alpha} y_{i} y_{j}\}.
			\end{eqnarray}	
 	Therefore, the determinant of an indefinite $(\alpha , \beta)$-metric is given by \cite{Sjors thesis}
	\begin{equation}\label{determinant}
\det (g_{ij})= \phi ^{n+1} \left(\phi -s  \phi' \right)^{n-2}  \left(\phi -s  \phi' + (sgn(\alpha ^2 )\,b^2 -s^2 )\phi'' \right)\det (\alpha_{ij}),
\end{equation}
where $b^2 := ||\beta||_{\alpha}^2 = \alpha ^{ij} b_i\, b_j$.
It should be noted that when $\alpha$ is a positive definite Riemannian metric, thus $sgn (\alpha ^2)=1$, the formulas \eqref{gen metric} and \eqref{determinant} reduces to the well known expressions of $(\alpha ,\beta )$-metric. In addition, the condition for the matrix $(g_{ij})$ to be positive definite matrix \cite[Lemma~1.1.2]{ChernShen2005}  is 
\begin{equation}\label{positive definite cond}
\phi (s)-s\phi'(s)+(b^2-s^2)\phi''(s)>0,\,\,\,\,\, |s| \leq \, ||\beta||_{\alpha}   . 
\end{equation}
	\subsection{Generalized $m$-Kropina metrics  }		
	\begin{Def} \label{Def:F}
	 Consider a smooth  function $\phi$  on its domain defined by
	 \begin{align}\label{eq:ppdsd}
         \phi(s)  =   s^{-m}( c+r s^2)^{\frac{1+m}{2}}, 
     \end{align}
     where $c,m,r$ are real constants and the $1$-form $\beta$ is nowhere vanishing.  The assumption $c+rs^2 > 0,\, s>0$  is necessary for any real exponent, appears in \eqref{eq:ppdsd},  to be defined (however,  this condition is not strictly necessary if $m$ is an integer). 
      The indefinite $(\alpha , \beta)$-metric  $F:= \alpha\,\phi(s)$  is called generalized $m$-Kropina metric  \emph{\cite{H1}}. More restrictions can be added, e.g.,  $mc \neq 0,\,\, rs^2 \neq cm ,\,\, m\neq 0,-1$ for certain reasons which will be clarified in the sequel.
\end{Def}  

Clearly, a generalized $m$-Kropina metric defined on the conic subbundle $$\mathcal{A}= \left\{ (x,y): \alpha (x,y) \neq 0,\,\, \beta (x,y) >0 ,\,\, c \alpha^{2}(x,y) + r \beta^{2}(x,y) >0, \,\, \det(g_{ij})(x,y) \neq 0 \right\}.$$
Moreover, it reduces to an m-Kropina (or, Bogoslovsky-Kropina) metric $F= \alpha \,  s^{-m}$, when $c=1,\, r=0$, in addition, if $m=1$, it reduces to a Kropina metric $F= \alpha \,  s^{-1}$. \\ 

We can explain the parameters $c,m,r$ appeared in \eqref{eq:ppdsd} as follows: 
\begin{itemize}
\item[(1)]
\textit{Anisotropy Degree $m$}: it quantifies Lorentz violation, this exponent dictates the fundamental severity of the local symmetry breaking.
As our rigidity results (appears in sequence) demonstrate,   the arithmetic nature of $m$,  fractional values heavily restrict the Einstein spacetime to Ricci flat. 
Also,  $m$ controls the nonlinear structure of the anisotropic contribution. \item[(2)]\textit{Coupling Weight $r$}: inspired by modified gravity theories like the SME, a particle's trajectory is perturbed by a background vector field (e.g., a Bumblebee vacuum expectation value, represented by the $1$-form $\beta$ \cite{Bluhm:2004ep, Jacobson:2000xp}).  It determines how strongly a test particle interacts with this preferred background direction $\beta$ relative to standard gravity.
\item[(3)]\textit{Isotropic Baseline $c$:} by the correspondence principle, any viable modified gravity theory must allow for a limit that recovers standard General Relativity. The parameter $c$ anchors the metric, ensuring that the $\alpha$ survives as the foundational gravitational baseline when local anisotropies are negligible.
\end{itemize}
\begin{Lem} \label{lemma:positive_definite_conic}
Let $F$ be a generalized $m$-Kropina metric 
with $m ,\,c > 0$. Then $F$ is a positive-definite conic Finsler metric with $0 < s \leq b$  if and only if 
\begin{equation}
    c + r b^2 > 0.
\end{equation}
\end{Lem}
\begin{proof}The condition for the matrix $(g_{ij})$  of an $(\alpha ,\beta )$-metric to be positive definite  \cite[Lemma~1.1.2]{ChernShen2005}  is given by  $\phi(s) > 0$ and \eqref{positive definite cond}.
Due to the $s^{-m}$ term, the generalized $m$-Kropina metric is singular at $s=0$. Thus, we restrict our analysis to the conic region $0 < s \leq b$. Because $c > 0$ and we will require $c + rb^2 > 0$, it follows that $c + rs^2 > 0$ for all valid $s$, ensuring $\phi(s) > 0$.
We have
\begin{equation}\label{phi'}
    \phi'(s) = s^{-m-1} (c+rs^2)^{\frac{m-1}{2}} (rs^2 - mc),
\end{equation}
\begin{align}
    \phi(s) - s\phi'(s) &= c(m+1) s^{-m} (c+rs^2)^{\frac{m-1}{2}},\label{eq:first_part}
\end{align}
and
\begin{equation} \label{phi''}
    \phi''(s) = c(m+1) s^{-m-2} (c+rs^2)^{\frac{m-3}{2}} (cm + rs^2). 
\end{equation}
Substituting \eqref{eq:first_part} and \eqref{phi''} into \eqref{positive definite cond}, we obtain:
\begin{equation} \label{eq:factored_inequality}
    c(m+1) s^{-m-2} (c+rs^2)^{\frac{m-3}{2}} \Big[ c m b^2 + \big[ c(1-m) + r b^2 \big] s^2 \Big] > 0.
\end{equation}
The sign of \eqref{eq:factored_inequality} depends  on  $f(s^2):=c m b^2 + \big[ c(1-m) + r b^2 \big] s^2$ which is linear with respect to $s^2$, its minimum on the interval $[0, b^2]$ must occur at one of the boundaries: As $s^2 \to 0^{+}$: $f(0) = cmb^2 > 0$ (since $c, m, b > 0$) and at $s^2 = b^2$: $f(b^2) =  b^2(c + rb^2)$.
Thus, we  have $b^2(c + rb^2) > 0$, which simplifies to $c + rb^2 > 0$. Since both boundary limits are strictly positive, $f(s^2) > 0$ for all $s^2 \in (0, b^2]$. 
\end{proof}

\begin{Lem} \label{lemma:non_degenerate_conic}
Let $F$ be a generalized $m$-Kropina metric, $n \geq 3$. Then $F$ is a non-degenerate (pseudo)-Finsler metric on the conic domain $0 < s \leq b$ if and only if the following conditions hold for all $s$ in the domain:
$$c \neq 0 ,\quad m \neq -1,\quad c + rs^2 > 0,\quad\, sgn(\alpha ^2 )\,c m b^2 + \big[ c(1-m)+ sgn(\alpha ^2 ) r b^2 \big] s^2 \neq 0.$$
\end{Lem}

\begin{proof}
For an $n$-dimensional $(\alpha, \beta)$-metric $F = \alpha \phi(s)$, the fundamental metric tensor $g_{ij}$ is non-degenerate if and only if its determinant (given by \eqref{determinant}) is non-zero. 
Thus, $F$ is non-degenerate if and only if  $$\phi(s) \neq 0, \quad \phi(s) - s\phi'(s)  \neq 0, \quad \phi (s)-s\phi'(s)+(sgn(\alpha ^2 )\,b^2-s^2)\phi''(s) \neq 0$$  on the domain  $0 < s \leq b$.  Thus, by  \eqref{eq:first_part} and \eqref{eq:factored_inequality} in case of $sgn(\alpha ^2 ) =1$, we get 
   \[c + rs^2 > 0, \quad c \neq 0,  \quad m \neq -1,\quad c m b^2 + \big[ c(1-m) + r b^2 \big] s^2 \neq 0.\]
In case of $sgn(\alpha ^2 ) =-1$,  the condition \[ \phi (s)-s\phi'(s)+(sgn(\alpha ^2 )\,b^2-s^2)\phi''(s) \neq 0\] reads
\[-c m b^2 + \left[ c(1-m) - r b^2 \right] s^2 \neq 0.\]
  \vspace*{-1.3cm}\[\qedhere\]
\end{proof}
\begin{Rem}\label{rem:spacetime conditions}\emph{\cite[\S 4.4]{Voicu F spacetime}}
To determine the spacetime conditions for the family of generalized $m$-Kropina metrics, taking the Lorentz signature of both $\alpha$ and $F$ to be $(+,-,-,-)$. 
The necessary and sufficient conditions for the generalized $m$-Kropina $4$-dimensional manifold  defines a Finsler spacetime (to describe a valid physical spacetime ensuring the causal cone of  future-pointing timelike vectors remains well-defined that prevent the causal structure from collapsing,  and non-degenerate Lorentzian signature be valid) are the coupling weight and isotopic baseline parameters must be restricted to $c > 0$, $r < 0$ along with $s^2 \in (0, -c/r),$ and one of the following situations occurs:
\begin{enumerate}   \item[(i)] $b^2 \ge 0$ and $m \in (-1, 1]$: in this case, the future-pointing time-like cones of $F$ are the connected components of the cones $c\, \alpha^2 +r\, \beta^2  > 0$ contained in the future-pointing time-like cones $\alpha$. \item[(ii)] $b^2 < 0$ and $m \in (-1, 0)$: in this case, the future-pointing time-like cones of $F$ are regions of the cones $c\, \alpha^2 +r\, \beta^2  > 0$ lying both in the future-pointing time-like cones $\alpha$ and in the half-space $\beta  > 0$.\end{enumerate}
\end{Rem}

\subsection{Kinematic Predictability and the AR-Finsler Property}
A significant mathematical challenge in modifying gravity via Finsler geometry is the highly non-linear, often non-polynomial nature of the field equations.
To ensure that the geometry yields physically meaningful and analytically stable equations of motion for test particles, we  analyze the rationality of the metric.
	\begin{Def}\emph{\cite{arFinslerTaha}}
		A conic pseudo-Finsler metric $F$ is said to be an almost rational Finsler (simply, AR-Finsler) metric if the components of its metric tensor  $g_{ij}(x,y)$ are  expressed in the form
		\begin{equation}\label{gdef}
			g_{ij}(x,y) =\eta(x,y)\, a_{ij}(x,y), \quad \forall (x,y) \in \mathcal{A}_{x}, \,\,\forall x \in M
		\end{equation} 
		where \begin{description}
			\item[(a)] $\eta: \mathcal{A}\longrightarrow \mathbb{R}$ is a smooth function,
			\item[(b)] the matrix $( a_{ij}(x,y))_{1 \leq i,j \leq n}$ is  symmetric non-degenerate with $a_{ij}(x,y)$ are rational functions in the fiber coordinate $y$ for all $i,j = 1,...,n$ (which means that $a_{ij}(x,y)$ can be written as a ratio of polynomials in the variable $y$ and no matter how $a_{ij}(x,y)$ looks like in the variable $x$).
\end{description}
		If, in addition,  $\eta$ is a rational function in $y$,   $F$ is called a rational Finsler metric.\footnote{It should be noted that $F \text{ rational in } y \Longrightarrow g_{ij} \text{ rational in } y,$ while the converse does not generally hold.}
	\end{Def}
      \begin{Lem}\label{Rem F not rat}
A generalized $m$-Kropina Finsler function $$F= \alpha\,\phi(s)=  \alpha s^{-m}( c+r s^2)^{\frac{1+m}{2}}$$ is a rational function in $y$ if $m$ is an odd integer. In particular, $F$ is an irrational function in $y$ provided that $m \in 2 \mathbb{Z}$. 
\end{Lem}
\begin{proof}
As $F$ can be written in the form
\begin{equation} \label{F interms of alpha , beta}
F=  \alpha s^{-m}( c+r s^2)^{\frac{1+m}{2}}=  \alpha \left(\frac{\alpha}{\beta}\right)^{m}\left( c+r \frac{\beta ^2}{\alpha ^2}\right)^{\frac{1+m}{2}} =  
   \frac{1}{\beta^{m}} \left( c \alpha ^{2}+r \beta ^2 \right)^{\frac{1+m}{2}}.
\end{equation} 
Since, $\beta$ is linear in $y$ and $\alpha ^2$ is quadratic in $y$, that is, both $\beta$ and $\alpha ^2$ are rational functions in $y$, then  $F$ is rational function in $y$ when $(m+1)$ is an even integer which means $m$  is an odd integer. 
\end{proof}
 \begin{Rem}
 It should be noted that if a Finsler structure  $F$ is a rational function in $y$, then $F^2$ is a rational function in $y$.  This is the case of generalized $m$-Kropina metric with $m$ odd integer. Further,  the partial differentiation $ \dot{\partial}_{i}$ of a geometric object that is rational in $y$ remains rational in $y$. Thereby,  the resulting Finsler metric tensor components $g_{ij}(x,y)$ are rational functions in $y$ and hence $F$ is a rational Finsler metric. However, the converse is false.  As an example, the $m$-th root metric is a rational Finsler metric when $m\geq 3$ and its Finsler function $F$ is not rational in $y$ \emph{\cite{arFinslerTaha}}. We show that $F$ with  $m \in 2\mathbb{Z}$ is a rational Finsler metric in the next result (Theorem \ref{generalized Kropina is AR }).   
 \end{Rem}
     \begin{Thm}\label{generalized Kropina is AR }
A generalized $m$-Kropina metric $F$ is an AR-Finsler metric for $m \in \mathbb{R}\setminus \mathbb{Z}$ and $F$ is a rational Finsler metric, $g_{ij}$ are rational functions in $y$,  for  $m \in \mathbb{Z}\setminus \{-1,0\}$.    
     \end{Thm}
     \begin{proof}
   According to \eqref{gen metric}, the Finsler metric tensor components of $F$  are given by 
	\begin{eqnarray} \label{gen m-Krop.}
	g_{ij}&=& \frac{c(m+1)( c+r s^2)^m}{s^{2m}}\, sgn (\alpha ^2)\,\alpha_{ij} \nonumber \\
	&& + \frac{(c^{2}m(2m+1) -cr(m-1)s^{2} +r^2 \, s^4)(c+rs^2 )^{m-1}}{s^{2m+2}} \, b_{i}\, b_{j} \\ \nonumber &&
	- \frac{2 c^{2}m(m+1)(c+rs^2 )^{m-1}}{\alpha \,s^{2m+1}} \, \left(sgn (\alpha ^2)\,\{b_{i} y_{j} + b_{j} y_{i} \} -  \frac{\beta}{\alpha ^2} y_{i} y_{j} \right). 
	\end{eqnarray}		
		Further,  \eqref{gen m-Krop.} can be written as  	
		\begin{eqnarray}	
		g_{ij}&=&\frac{(c+rs^2 )^m}{s^{2m}}\,\{ c(m+1)\,sgn (\alpha ^2)\, \alpha_{ij} + \frac{(c^{2}m(2m+1) -cr(m-1)s^{2} +r^2 \, s^4)}{s^2\,(c+rs^2 )}  \, b_{i}\, b_{j} \nonumber\\  &&
		-\frac{2 c^{2}m(m+1)}{\beta (c+rs^2 )} \, \left( b_{i} y_{j} + b_{j} y_{i}  -  \frac{\beta}{\alpha ^2} y_{i} y_{j} \right)\}.
		\end{eqnarray}
		Thereby, expression of $g_{ij}$ given by \eqref{gen m-Krop.} has the form  \eqref{gdef} with
	\begin{eqnarray}
	 \eta &= &\frac{(c+rs^2 )^m}{s^{2m}}= \left(\frac{c\alpha^2 +r\beta{^2} }{\beta{^2}}\right)^m = (c+rs^2 )^{-1}\, \phi^{2} (s); \label{eta}\\
	 a_{ij}&=& c(m+1)\,sgn (\alpha ^2)\, \alpha_{ij} + \frac{(c^{2}m(2m+1) -cr(m-1)s^{2} +r^2 \, s^4)}{s^2\,(c+rs^2 )}  \, b_{i}\, b_{j}\nonumber \\  &&
		-\frac{2 c^{2}m(m+1)}{\beta (c+rs^2 )} \, \left(sgn (\alpha ^2)\,(b_{i} y_{j} + b_{j} y_{i} ) -  \frac{\beta}{\alpha ^2} y_{i} y_{j} \right). \label{a_ij}
\end{eqnarray}		
Since both $b_{j}$ and $\alpha_{ij}$ are functions of $x$ only, $s^{2} = \frac{\beta ^2}{\alpha ^2}$ is rational functions in $y$, the functions $a_{ij}(x,y)$ are rational in $y$ for all real values of $m$. Moreover, when $m \in \mathbb{Z}\setminus \{-1,0\}$, the function $\eta$ is  rational in $y$, however,  when $m \in \mathbb{R} \setminus \mathbb{Z}$, the function $\eta$ is  irrational in $y$. Hence, a generalized $m$-Kropina metric is an AR-Finsler metric. 
     \end{proof}
     
     \begin{Rem}\label{Cartan tensor}\footnote{All over the paper we exclude $c=0,\, m=-1$.}
     $(1)$ From \eqref{phi'}, $\phi '(s) =0 \Longleftrightarrow rs^2 -cm =0.$ Therefore, the condition $rs^2 \neq cm$ is necessarily to make $\phi$ not constant or to make $F$ not a pseudo-Riemannian metric.\\
     $(2)$ The real constant $c=0$, from \eqref{F interms of alpha , beta},  implies $F=r \beta$ which means $F$ is linear and hence degenerate. Thus,  $c = 0$ has to be excluded. \\
     Similarly, $m \neq -1$ is required as $m=-1$ gives $F=  \beta$. \\
     $(3)$  The real constant $r=0$, gives $F= c^{\frac{m+1}{2}} \beta ^{-m} \alpha ^{m+1} = c^{\frac{m+1}{2}} \alpha s^{-m}$ which is a constant multiple of the $m$-Kropina metric.\\
     $(4)$ The real constant $m=0$,  from \eqref{eq:ppdsd},  gives $F=  \alpha \sqrt{c+rs^2}$ and $\eta =1$\emph{ (}by \eqref{eta}\emph{)} along with $a_{ij} =  sgn(\alpha ^2 ) \,c\,\alpha _{ij} + r b_i \, b_j $.  Consequently,  $g_{ij} = a_{ij}$ are functions of $x$ only and hence,  this $F$ is pseudo-Riemannian metric.
     \end{Rem}
 

     \subsection{Rationality of the associated Finslerian  objects}
     \begin{Lem}\label{main Cartan torsion}
    The Cartan tensor components $C_{ijk}$ of  a generalized $m$-Kropina metric $F$ are rational  in $y$ for $m \in \mathbb{Z}$. 
    However, its mean Cartan tensor  $I$ is a rational  function in $y$ for all real values of $m$.  \end{Lem}
    \begin{proof}
    Since,  a generalized $m$-Kropina metric $F$ is an AR-Finsler metric, then $C_{ijk}$ are given by \[C_{ijk}= \frac{1}{2} (a_{ij}\, \dot{\partial}_{k} \eta + \eta \,\dot{\partial}_{k}a_{ij}).\]
    As $a_{ij}(x,y)$ are rational in $y$ for all real values of $m$ by \eqref{a_ij} and $\eta$ is  a rational function in $y$ when $m \in \mathbb{Z}$ by \eqref{eta}.  Hence,  $C_{ijk}$ are rational  in $y$ for $m \in \mathbb{Z}$.\\
 The metric $g_{ij} =\eta \, a_{ij}$ has inverse $g^{jk}= \frac{a^{jk}}{\eta},\,\, \eta \neq 0 $ with $a^{jk}$ is the inverse of $a_{ij}$,
  thus we get
  \[I_{l}=g^{jk}\,C_{ljk}= \frac{a^{jk}}{\eta} \left( \eta\, \dot{\partial}_l a_{jk} + a_{jk}\, \dot{\partial}_l \eta \right)= a^{jk}  \dot{\partial}_l a_{jk} + n \, \dot{\partial}_l \log(\eta) .
		\]
		By direct calculations, we obtain  \[\dot{\partial}_{i} \log(\eta)= -\frac{2mc \left(\alpha^2 \, b_{i}- \beta \, \alpha_{ij}\, y^{j} \right)}{(c\alpha^2+r \beta ^2)\beta}.\] Therefore, functions $\dot{\partial}_{i} \log(\eta)$ are rational functions in $y.$\footnote{this result holds for positive definite  AR-Finsler metrics, see e.g. \cite[Proposition 3.4]{arFinslerTaha}} Hence, the mean Cartan torsion has components $I_{l}$ which are  rational  in $y$.
    \end{proof}
     \begin{Lem} \label{hij} Let $F$ be a generalized $m$-Kropina metric. Then, the angular metric tensor components $h_{ij}$ are rational functions in $y$ for  $m \in \mathbb{Z}$.
\end{Lem} 
\begin{proof} 
Assume that $m\in \mathbb{R}$.  Since, \[\dot{\partial}_{i} s= b_i \alpha^{-1} - \beta \alpha ^{-3} \alpha_{ij}\, y^{j}, \qquad \phi '(s)=\phi (s)\,\frac{r s^2 -c m}{c s+r s^3}.\]
  The associated normalized supporting element $\ell _{i} (x,y)$ are given by
  \begin{eqnarray*}
  \ell _{i} &=& \dot{\partial}_{i} F = \dot{\partial}_{i} (\alpha \, \phi (s))= (\dot{\partial}_{i} \alpha )\, \phi (s) + \alpha \, \phi '(s)\,\dot{\partial}_{i}s\\
  &=&  \phi (s)\, \left(\dot{\partial}_{i} \alpha  + \alpha \,\frac{r s^2-c m}{c s+r s^3} \,\dot{\partial}_{i}s \right)\\
  &=& 
\alpha^{-1} \,  \phi (s)\, \left(sgn(\alpha ^2 )\, \alpha_{ij}\, y^{j}+ \alpha \,\frac{r s^2-c m}{c s+r s^3} \,(b_i - \beta \alpha ^{-2} \alpha_{ij}\, y^{j})  \right).
\end{eqnarray*}
As $\phi (s)$ can be written in the form
 $$\phi (s)= \alpha^{-1} \,F= \alpha^{-1} \, \beta^{-m}\left( c \alpha ^{2}+r \beta ^2 \right)^{\frac{1+m}{2}}$$ and $$\frac{r s^2-c m}{c s+r s^3}= \alpha \,\frac{(r\beta ^2 -cm \alpha^2)}{c\beta \alpha ^2 + r \beta ^3},$$ we get
 \begin{eqnarray*}
  \ell _{i}  &=& \alpha^{-2} \, \beta^{-m}\left( c \alpha ^{2}+r \beta ^2 \right)^{\frac{1+m}{2}}\left( \alpha_{ij}\, y^{j}+ \alpha^2 \,\left(\frac{r\beta ^2 -cm \alpha^2}{c\beta \alpha ^2 + r \beta ^3}\right) \,(b_i - \beta \alpha ^{-2} \alpha_{ij}\, y^{j})  \right).
   \end{eqnarray*}
   Therefore,  $ \ell _{i}$ is rational in $y$  when $(m+1)$ is an even integer, or equivalently $m$  is an odd integer.  Consequently, the geometric quantity $\ell _{i} \ell _{j}$  is rational in $y$ when $m \in \mathbb{Z}$.  Hence, $h_{ij} = g_{ij} - \ell _{i} \ell _{j} $ are rational functions in $y$ (for $m \in \mathbb{Z}$). 
  \end{proof}   
 
 \begin{Thm} \label{rational geometric objects}
The geodesic spray  coefficients of a generalized $m$-Kropina metric are rational in $y$.  Moreover, each of the following Finslerian geometric objects associated with a generalized $m$-Kropina metric is rational in $y$: Barthel connection coefficients $N^{i}_{j}$, Berwald connection coefficients $G^{i}_{jk}$, Berwald curvature components $G^{i}_{jkl}$, mean Berwald curvature components $E_{ij}$, the Barthel connection curvature components
$R^{i}_{jk}$, Riemannian curvature components $R^{i}_{j}$,  Finsler Ricci tensor components $R_{ij} $and Finsler Ricci scalar $Ric $.
\end{Thm}    
\begin{proof}First, we prove that the spray coefficients $G^{i} $ of  a generalized $m$-Kropina metric $F$ are rational functions in $y$.
Since $F= \frac{1}{\beta^{m}} \left( c \alpha ^{2}+r \beta ^2 \right)^{\frac{1+m}{2}}$ is an $(\alpha ,\beta)$-metric,  $G^{i} $ can be calculated from the expression \cite{ChernShen2005}
\begin{equation}\label{spray of alpha beta metric}
G^{i} = G^{i}_{\alpha} + \alpha \, Q \, s^{i}_{0} + (r_{00}-2 \alpha \, Q \, s_{0} )( \Psi \, b^{i} + \alpha ^{-1} \, \Theta \, y^{i}),
\end{equation}
where 
$$
r_{00}:=\frac{1}{2}(b_{i;j}+b_{j;i})\,y^iy^j,\quad s^{i}_0:=\frac{1}{2}(b_{j;k}-b_{k;j})\,\alpha^{ij}\,y^k,\quad s_0:=y^i b^j \frac{1}{2}(b_{j;i}-b_{i;j}), $$
\begin{eqnarray*}
Q(s)&:= & \frac{\phi '(s)}{\phi (s)-s \phi '(s)} =\frac{r s^2 -c m}{c (m+1)s} \\
 \Theta (s) &:= & \frac{\phi (s) \phi '(s)-s \left(\phi (s) \phi ''(s)+\phi '(s)^2\right)}{2 \phi (s) \left[\left(b^2-s^2\right) \phi ''(s)+\phi (s)-s \phi '(s)\right]} = \frac{c m s}{c (m-1) s^2-b^2 \left(c m+r s^2\right)}\\
 \Psi (s) &:= & \frac{\phi ''(s)}{2 \left(\left(b^2-s^2\right) \phi ''(s)+\phi (s)-s \phi '(s)\right)} =\frac{c m+r s^2}{2 s^2 \left(b^2 r-c( m-1) \right)+2 b^2 c m}
\end{eqnarray*}
Therefore,
\small{
\begin{equation}\label{spray of gen m-Kropina}
G^{i} = G^{i}_{\alpha} + \alpha\frac{r s^2 -c m}{c (m+1)s} \, s^{i}_{0} + (r_{00}-2 \alpha   \frac{r s^2 -c m}{c (m+1) \, \beta}\, s_{0} )( \Psi \, b^{i} +  \frac{c m \beta}{c (m-1) s^2-b^2 \left(c m+r s^2\right)} \, y^{i}),
\end{equation}}
As $b_{i;j}$ denotes the covariant derivative of $b_i$ with respect to $\alpha$, then $b_{i;j}$ are functions on $x$ only. Also, as $ G^{i}_{\alpha} $ are  the geodesic spray coefficients of the Riemannian metric $\alpha$, they are quadratic in $y$ as well as $r_{00}$. Thereby, each of the functions $ s^{i}_0 ,\,\,  s_0 $ are linear in $y$. Similarly,  \[\alpha \, Q = \frac{r \beta^2 -c m  \alpha^2 }{c (m+1) \, 
\beta}, \qquad\alpha ^{-1} \Theta = \frac{c m \beta}{c (m-1) s^2-b^2 \left(c m+r s^2\right)}\]
are rational functions in $y$. In addition, the product of rational functions results a rational function as well as the sum of rational functions gives  a rational function.  Hence, from \eqref{spray of gen m-Kropina},  $G^{i} $ are rational functions in $y$.
\par The rest of the proof follows from the facts that the partial derivative with respect to $y$ of a geometric object which is rational in $y$ remains rational in $y$ (for example,   $N^{i}_{j}= \dot{\partial}_{j} G^i, \,\, G^{i}_{jk} =\dot{\partial}_{k} N^{i}_{j},\,G^{i}_{jkl} = \dot{\partial}_{l} G^{i}_{jk} \text{ and } E_{ij} = \frac{1}{2} G^{k}_{kij} $ are rational functions),  the partial derivatives $\dot{\partial}_i$ and  $\partial_{l}$ commute and the partial derivative with respect to $x$ of any rational function in $y$  remains rational in $y$ along with the multiplication or sum of rational functions always gives a rational function.  Thus,  from \eqref{curvature of Barthel connection}, the Barthel connection curvature components
 are rational functions in $y$. Therefore,  from \eqref{h Riem curv}  and \eqref{Ricci scalar and tensor} the functions $ R^i_k , \, R_{ij},\, Ric$ are rational in $y$.
\end{proof}
\begin{Rem}
An example of Finsler metric which have rational spray coefficients  but they are  not  AR-Finsler metrics is Randers metrics of Berwald type (which  has rational  spray coefficients)\emph{ \cite{arFinslerTaha}.}
\end{Rem}

\subsection{Rigidity results and special Finsler metrics}
Here, we make use of the rationality of the associated Finslerian geometric objects to get some rigidity results and classify the family of  generalized $m$-Kropina metrics.
\begin{Thm}\label{main results Einstein}
If $F$ is a  generalized $m$-Kropina metric,  with $m \notin \mathbb{Z}$,  of an Einstein type, then $F$ is Ricci-flat. 
	\end{Thm}
	\begin{proof}
	 Suppose that a generalized $m$-Kropina metric $F$ is an Einstein metric, that is,  for some $K \in C^{\infty}(M)$, we have
	  $$Ric(x,y) = (n-1)  K(x) \,F^2 (x,y).$$
Now, $F^2 $ is irrational in $y$ for $m \notin \mathbb{Z}$ as  the function $\eta$, given by \eqref{eta}, is.  However,  $Ric$ is a rational function in $y$ by Theorem \ref{rational geometric objects}. That is, the L.H.S.  is rational function in $y$ while the R.H.S is irrational in $y$. Hence, each of $Ric(x,y) $ and $K(x)$ vanishes identically which means $F$ is Ricci-flat.
	\end{proof}
	It should be noted that Theorem \ref{main results Einstein} holds for AR-Finsler metrics with $\partial_{i} \log (\eta(x,y))$ being rational functions in $y$,   see e.g. \cite[Theorem 3.5]{arFinslerTaha}.
	\begin{Thm}\label{main results}
	Let $F$ be a generalized $m$-Kropina metric with $m \in 2\mathbb{Z}$.  We have the following:
	\begin{description}
\item[(a)] If $F$ is a weakly Einstein metric,  then $F$ is an Einstein metric. 	
\item[(b)]  If  $F$ has isotropic mean Landsberg curvature,  then $F$ is  weakly Landsbergian. 
\item[(c)]  If  $M$ equipped with an arbitrary volume form $dV_\mu$ and $F$ has isotropic $S_{\mu}$-curvature,  then the $S_{\mu}$-curvature vanishes identically. 
\item[(d)] If $F$ has isotropic mean Berwald curvature,  then $F$ is  weakly Berwaldian.
\item[(e)]  If $F$ has relatively isotropic Landsberg curvature,  then $F$ is Landsbergian.
\item[(f)]  If $F$ has almost vanishing $\mathbf{H}$-curvature,  then $\mathbf{H}=0$.
\item[(g)]  If $F$ has almost isotropic flag curvature $\mathbf{K}$ provided that $n \geq 3$, then $\mathbf{K}$ is constant. 
\item[(h)] \label{Minkowskian}If $n\geq 3$ and $F$ has scalar flag curvature $\mathbf{K}$ and isotropic $S_{\mu}$-curvature, then the $S_{\mu}$-curvature vanishes identically and $F$ is locally Minkowskian. 
\end{description}
\end{Thm}
	\begin{proof}
	By Lemma \ref{Rem F not rat},  a generalized $m$-Kropina metric $F$ is an irrational function in $y$ when $m \in 2 \mathbb{Z}$.  However,  $F^2 $ is rational in $y$ for $m \in \mathbb{Z}$ (as the function $\eta$ given by \eqref{eta} is), this is because $F^2 = g_{ij} y^i y^j = \eta a_{ij} y^i y^j $.
	\begin{description}
	\item[(a)] Now,  assume that $F$ is  weakly Einstein, that is,  $$Ric = (n-1) \left\{ \frac{3 \theta}{F} + K \right \}F^2,$$
for some $1$-form $\theta$ on $M$ and for some smooth function $K$ 	on $M$. 
	  Which is equivalent to
	\begin{equation}\label{weak Eins}
Ric(x,y)  - (n-1) K(x) F^{2}(x,y)= 3(n-1)  \theta_{l} (x)\, y^l \, F(x,y).
	\end{equation}
	 As $Ric$ is a rational function in $y$ by Theorem \ref{rational geometric objects}, the L.H.S.  is a rational function in $y$ while the R.H.S is an irrational function in $y$. Therefore, for  $F\neq 0$, the $1$-form $\theta$ vanishes identically and hence $F$ is an Einstein metric.
\item[(b)]	Since,  the mean Cartan torsion  $I(x,y)$ is rational  in $y$ for all real values of the anisotropy degree $m$,  by Lemma  \ref{main Cartan torsion}.  Moreover, $$J :=J_k\, dx^k ,\,\quad J_{k}= g^{ij}\,L_{ijk} = y^p\, \partial_p I_k -2 \,G^p\, \dot{\partial}_{p} I_k - N^{p}_{k}\, I_p $$
is rational function in $y$ as $G^p ,\, N^{p}_{k}$ are rational functions in $y$ (by Theorem \ref{rational geometric objects}).   Assume that $F$ has isotropic mean Landsberg curvature, i.e., for some $A \in C^{\infty}(M)$ we have,  $$J(x,y) = A(x)\, F(x,y) \, I(x,y).$$
 Hence,  both sides of the last equation vanishes  and consequently $J=0$.
\item[(c)] Assume that $F$ has an isotropic $S_{\mu}$-curvature, that is,  for some $A \in C^{\infty}(M)$,
$$S_{\mu}(x,y)= (n+1)  F(x,y)\, A(x).$$ It is clear that $S_{\mu}(x,y)$ is a rational function in $y$, as $$S_{\mu}(x,y)=N^i_i (x, y)-y^i\, \partial_{i}\log\left( \sigma_{\mu}(x)\right).$$ Therefore,  for $F(x,y)  \neq 0$, the function $A$  vanishes identically, or equivalently, $S_{\mu}$-curvature vanishes.
\item[(d)] Suppose $F$ has isotropic mean Berwald curvature, that is, $$FE_{ij} = \frac{n+1}{2}f\, h_{ij},$$ for some $ f \in C^{\infty}(M)$. By Theorem \ref{rational geometric objects}, $E_{ij}$ are rational functions in $y$. Also, $h_{ij}$ are rational functions in $y$ for all $m \in \mathbb{Z}$ (by Lemma \ref{hij}). However, $F$ is an irrational function in $y$ for $m \in 2 \mathbb{Z}$, by Lemma \ref{Rem F not rat}. Hence, we get $E_{ij} =0$ and $f=0$ which means $F$ is  weakly Berwaldian.
\item[(e)]  Suppose $F$ has relatively isotropic Landsberg curvature, that is, $$L_{ijk} = F\,f\, C_{ijk},$$for some $f \in C^{\infty}(M)$.
Since, $C_{ijk}$ are rational functions in $y$ for all $m \in 2\mathbb{Z}$ from Lemma \ref{main Cartan torsion}. Thus, $F\,f\, C_{ijk}$ are  irrational functions  in $y$ for $m \in 2\mathbb{Z}$. On the other hand, the components of the Landsberg curvature $$L_{ijk}=-\frac{1}{4} y^{q} g_{ql}G^{l}_{ijk} = -\frac{1}{4} y^{q} \eta \,a_{ql}G^{l}_{ijk}$$ are rational in $y$, for $m \in \mathbb{Z}$  (as $\eta$ is rational in $y$ and $G^{l}_{ijk}$ are rational in $y$ by Theorem \ref{rational geometric objects}). Hence,  $L_{ijk}=0$ which means $F$ is Landsbergian.
\item[(f)] Assume that  $F$ has almost vanishing $\mathbf{H}$-curvature, that is, $$FH_{ij} = \frac{(n+1)}{2} A\, h_{ij}, \quad \text{  for some } A \in C^{\infty}(M).$$We have $h_{ij}$ are rational functions in $y$ for all $m \in \mathbb{Z}$  (by Lemma \ref{hij}).  Thus, $\frac{(n+1)}{2} A(x) h_{ij}(x,y)$ are rational functions in $y$ for $m \in \mathbb{Z}$.  Also, for $m \in 2\mathbb{Z}$, the function $F$ is  irrational in $y$. 
\par However, $H_{ij}$ are rational in $y$ as $$H_{ij} = y^l ( \delta _{l} E_{ij} - G^{k}_{il} E_{kj}- G^{k}_{jl} E_{ki})=y^l ( \partial_{l} E_{ij} - N^{k}_{l} \,\dot{\partial}_{k} E_{ij} - G^{k}_{il} E_{kj}- G^{k}_{jl} E_{ki}) $$  by making use of Theorem \ref{rational geometric objects}.  Thereby,  the functions $FH_{ij}$ are irrational in $y$. Hence,  we obtain $H_{ij}=0$ and $A=0$ (as $h_{ij}$ can not be vanish identically) for $m \in 2\mathbb{Z}$. Consequently, $\mathbf{H}=0$ for $m \in 2\mathbb{Z}$. 
\item[(g)] Assume that $F$ has almost isotropic flag curvature $\mathbf{K}$, that is, $$\mathbf{K}(x,y)= \frac{y^i \partial_{i} A(x)}{F(x,y)} + \sigma (x) .$$
The flag curvature $\mathbf{K}= \mathbf{K}(P,y)$ of a flag $(P,y)$, where $P=span\{ y,u\} \subset T_{x}M$ is defined by 
\begin{equation}\label{flag curvature def}
\mathbf{K}= \frac{g_{ij} \, R^{i}_{k}(x,y)\, u^j u^k}{F^{2}(x,y)\, g_{ij} \,u^i u^j -[g_{ij}(x,y) \,y^i u^j]^2 }.
\end{equation}
As $F$ is a rational Finsler metric for $m \in \mathbb{Z}$ (by Theorem \ref{generalized Kropina is AR }),  that is,  $g_{ij} =\eta a_{ij} $ and $F^2 = \eta \,a_{kl}\,y^k y^l $, thus we get
\begin{eqnarray*}
\mathbf{K}&=& \frac{\eta a_{ij} \, R^{i}_{k}(x,y)\, u^j u^k}{\eta a_{kl}\,y^k y^l \, \eta a_{ij} \,u^i u^j -[\eta \, a_{ij}(x,y) \,y^i u^j]^2 }=\frac{1}{\eta} \frac{ a_{ij} \, R^{i}_{k}(x,y)\, u^j u^k}{ a_{kl}\,y^k y^l \, a_{ij} \,u^i u^j -[ a_{ij}(x,y) \,y^i u^j]^2 }
\end{eqnarray*}
Since, $R^{i}_{k}(x,y)$ are rational functions on $y$ from Theorem \ref{rational geometric objects} and $\eta$ is a rational function in $y$ when $m \in \mathbb{Z}$ from \eqref{eta}, we deduce that $\mathbf{K}$ is rational function in $y$ for all  $m \in \mathbb{Z}$. \\ Thus, the function $\mathbf{K} -\sigma $ is still rational in $y$. However, for $m$ even integer,  $\frac{y^i \partial_{i} A(x)}{F}$ is an irrational function in $y$. Hence, $y^i \partial_{i} A(x) = 0$. Consequently, $\mathbf{K}= \sigma$. Therefore, $\mathbf{K}$ is constant (by Schur's Lemma as $n \geq 3$).
\item[(h)] Assume that $F$ has isotropic $S_{\mu}$-curvature. Then,  item \textbf{(c)} of this theorem, $S_{\mu}$-curvature vanishes identically. Also, assume that $F$ has  scalar flag curvature. \\
 Since a generalized $m$-Kropina metric $F$ is an $(\alpha , \beta)$-metric not of Randers type, we can apply the main result of \cite{Cheng2010} (which holds for pseudo-Finsler metrics) \lq \lq  any $(\alpha , \beta)$-metric not of Randers type has scalar flag curvature and vanishing $S_{\mu}$-curvature if and only if it is Berwaldian with vanishing flag curvature." and deduce that   $F$ is Berwaldian with vanishing flag curvature. Hence, $F$ is  locally Minkowskian.
\vspace*{-0.6cm}\[\qedhere\]
\end{description}
\end{proof}
It should be noted that Theorem \ref{main results} holds for rational Finsler metrics (i.e,  $\eta$ rational function in $y$) with irrational Finsler function $F$ in $y$ and  $\partial_{i} \log (\eta(x,y))$ being rational functions in $y$.
\begin{Rem}
The rationality arguments in this section are independent of the signature, provided the relevant nondegeneracy and conic domain assumptions are satisfied, that is, the results obtained are true for positive definite (non-degenerate, with certain signature, etc)  generalized $m$-Kropina metrics as their proofs depend mainly on the rationality of the geometric objects.  
\end{Rem}

Let us conclude this section with a summary of  the rational properties of Finslerian geometric objects associated with a generalized $m$-Kropina metric.

\begin{center}
\begin{center}
\textbf{Table: Rationality of the objects of a generalized $m$-Kropina metric\footnote{All statements are understood on the nondegenerate conic domain and with the exceptional values excluded as specified in Definition \ref{Def:F} and Remark \ref{Cartan tensor}.}}
\end{center}
\vspace{10pt}
\begin{tabular}{|c|c|} \hline
\textbf{Geometric object} & \textbf{Rationality in $y$ under a condition on $m$} \\ \hline
$F,\,l_{i}$  & m: odd integer \\ \hline
$\eta,\, g_{ij},  C_{ijk},\,h_{ij}$  & $m \in \mathbb{Z}$\\  \hline
$\partial_{j} \log(\eta),\,\dot{\partial}_{j} \log(\eta)$   & $m \in \mathbb{R}$  \\  \hline
$I_i ,\,G^i ,\, N^i_j ,\, G^i_{jk},\, G^i_{jkl}$ &  $m \in \mathbb{R}$   \\  \hline
$E_{ij},\, H_{ij},\, J_i $ &  $m \in \mathbb{R}$   \\ \hline
$R^i_j,\, Ric,\, S_{\mu}$-curvature  &  $m \in \mathbb{R}$   \\ \hline
$\mathbf{K},\,L^i_{jkl}$  & $m \in \mathbb{Z}$ \\  \hline
\end{tabular}
\end{center}

\section{Application  to Finsler gravity}
There are several attempts to derive Einstein's field equations in the context of Finsler geometry. We consider only one of them, then using the results mentioned above we determine under what conditions a generalized $m$-Kropina metric $F$ becomes an exact solution to it.

\par Finslerian version of Einstein's field equation (in vacuum, i.e.,  without cosmological constant) is derived by Pfeifer-Wohlfarth in \cite{Pfeifer:2011xi,Hohmann_2019} and proved to be a good candidate for describing  Finsler gravity see, e.g. \cite{Sjors thesis, ref5}. For, $n=4$, it is given by
\begin{equation}\label{general field eq}
-2 Ric +\frac{2F^2 }{3}g^{ij} R_{ij}+\frac{2F^2 }{3}g^{ij} \{\dot{\partial}_{i}(y^l \delta_l J_j -N^l _j J_l) +J_{j |i}\}=0,
\end{equation}
where  \lq \lq $|$" is the horizontal covariant derivative with respect to Chern connection.

\begin{Rem} $(1)$ One can show that, when $F$ is a pseudo-Riemannian metric:\\
 The field equation \eqref{general field eq} reduces to Einstein's field equation in vacuum, $R_{ij}=0$. \\
$(2)$ If $\alpha$ is a solution to Einstein's field equations in vacuum (i.e. $\alpha$ is Ricci-flat) and the $1$-form $\beta$ is covariantly constant with respect to $\alpha$ (that is, $b_{i;j} =0$) then the constructed $(\alpha , \beta)$-metric $F$ is an exact solution to \eqref{general field eq}. This can be shown by plugging $b_{i;j} =0$ into \eqref{spray of alpha beta metric},  yields  $r_{00}= 0,\,\, s^i_0 =0,\,\,s_{0}=0$. Thereby, $G^i = G^i_{\alpha}$. Hence, $F$ is Berwaldian (which implies $J_k =0$) and consequently, the associated Ricci tensors coincide which gives $F$ is Ricci flat as $\alpha$ is an exact solution of Einstein's field equations in vacuum \cite{Sjors thesis}. \\
\end{Rem}
\begin{Pro}
 \begin{description}
\item[(a)] For a Finsler metric of weakly Einstein type,  the field equation  \eqref{general field eq} has the form
 \begin{equation}
2KF^2  -3\theta \,F+\frac{2F^2 }{3}\, g^{ij} \{\dot{\partial}_{i}(y^l \delta_l J_j -N^l _j J_l) +J_{j |i}\} = 0.
\end{equation}
\item[(b)]For a Finsler metric of Einstein type,  the field equation  \eqref{general field eq} has the form
 \begin{equation}
3K \, F^2 +F^2\, g^{ij} \{\dot{\partial}_{i}(y^l \delta_l J_j -N^l _j J_l) +J_{j |i}\} = 0.
\end{equation}
\end{description}
\end{Pro}
\begin{proof}\begin{description}
\item[(a)] Let $F$ be  a weakly Einstein Finsler metric, that is,  $$Ric = (n-1) \{ 3 \theta \,F + K \,F^2 \},$$
for some $1$-form $\theta$ on $M$ and $K \in C^{\infty}(M)$.  Thereby,  by \eqref{Ricci scalar and tensor}, we get
  \[R_{ij}= (n-1)\left\{\frac{3}{2} \left(\theta_i \ell_j  +\theta_j \ell_i+ \theta \, \frac{h_{ij}}{F}\right)+  K \, g_{ij}\right\}.\] Therefore, as $\ell_i = \frac{y^i}{F}$ and for $n=4$, we have $g^{ij}h_{ij} =3,\, g^{ij}g_{ij}=4$:
\small{  \begin{eqnarray*}
  -2 Ric +\frac{2F^2 }{3}g^{ij} R_{ij} &=& -6\{ 3 \theta \,F + K \,F^2 \} +\frac{2F^2 }{3}\, g^{ij} \{\frac{9}{2}(\theta_i \ell_j  +\theta_j \ell_i+ \theta \, \frac{h_{ij}}{F})+ 3 K \, g_{ij} \} \\ &=& -18\theta \,F - 6K \,F^2  + F^2  3(2 \theta_i \ell^i + \theta \, \frac{h_{ij} g^{ij}}{F})+ 8K F^2  
   \\ &=& 2KF^2  -18\theta \,F +6\theta \,F +9\theta \,F =2KF^2  -3\theta \,F
  \end{eqnarray*}}
\item[(b)] Let $F$ be a Finsler metric of Einstein type, that is for some $K \in C^{\infty}(M)$, we have $ Ric = (n-1)  K \, F^2 $. Thus,  by \eqref{Ricci scalar and tensor}  $R_{ij}= (n-1)  K \, g_{ij} $. Therefore, for $n=4$,
\[-2 Ric +\frac{2F^2 }{3}g^{ij} R_{ij} = -6K \, F^2  + 2F^2 g^{ij}\, K \, g_{ij}  = -6K \, F^2 +8K \, F^2 =2K \, F^2 .\]
\end{description}
The proof is completed by substituting  in   \eqref{general field eq}.
\end{proof}
The following theorem gives conditions under which generalized $m$-Kropina metrics, which are not necessarily Berwaldian, provide exact solutions to the Finslerian version of Einstein's field equations \eqref{general field eq}.
\begin{Rem} Every generalized $m$-Kropina metric $F$ of Einstein type with fractional anisotropic degree $m$ is an exact solution to the field equation \eqref{general field eq}, provided that the spacetime conditions in Remark \emph{\ref{rem:spacetime conditions}} are applied,  if $F$ satisfies
\begin{equation}\label{Ricci flat feild eq}
F^2 \, g^{ij} (2 J_{i|j} + y^l \delta_{l} \dot{\partial}_{i}J_{j} + N^{k}_{i}\dot{\partial}_{k} J_{j}-N^{k}_{j} \dot{\partial}_{i}J_{k})=0.
\end{equation}
This can be seen from Theorem \ref{main results Einstein},   $F$ is Ricci-flat (that is, $Ric=0$ which is equivalent to $R_{ij}=0$).  Thereby, if $F$ satisfies \eqref{Ricci flat feild eq}, then it is an exact solution to \eqref{general field eq}. This  isolates the necessary and sufficient residual condition for an Einstein generalized $m$-Kropina metric with  fractional $m$ to solve the Pfeifer–Wohlfarth equation. It should be noted that any Finsler metric of weakly Landsbergian type are examples of  exact solutions of \eqref{Ricci flat feild eq}.
\end{Rem}
\subsection{Phenomenological Case Studies}
To demonstrate the direct utility of generalized $m$-Kropina metrics in modern modified gravity, we map the abstract metric parameters $(m, r, c)$ to established phenomenological constants. 

\subsubsection{The pseudo-Riemannian and Disformal Gravity Limit ($m=0$)}
To satisfy the correspondence principle, a viable modified gravity theory must smoothly recover standard quadratic kinematics in the limit where Finslerian Lorentz violation vanishes. In our generalized $m$-Kropina  metric, this boundary occurs exactly at $m=0$. In attempts to model galactic dynamics without dark matter, Tensor-Vector-Scalar (TeVeS) gravity relies heavily on Bekenstein's disformal transformations \cite{Bekenstein:2004ne}. In these models, matter couples to a physical metric constructed from a pseudo-Riemannian base $\alpha$ and a dynamical vector field $U_\mu$. 

Our $m$-Kropina framework natively embeds this disformal gravity as its zero anisotropic degree limit ($m=0$), the Finsler metric tensor $g_{ij}$ exactly reduces to:
\begin{equation*}
g_{ij} = c \, sgn(\alpha^2) \alpha_{ij} + r b_i b_j
\end{equation*}
At this exact limit, all higher-order directional dependence vanishes.  The spacetime ceases to be a true Finsler manifold has the mathematical form of a disformal pseudo-Riemannian metric of Bekenstein type, where $c$ and $r$ act as the conformal and disformal coupling constants. Thus, $m=0$ represents the exact mathematical boundary where Lorentz violation vanishes and standard, albeit disformally modified, GR kinematics resume.
\subsubsection{Spontaneous Lorentz Breaking: Bumblebee and Einstein-Aether Gravity}
In the Standard-Model Extension and Einstein-Aether theory, Lorentz symmetry is spontaneously broken by a vector field (such as the Bumblebee field $B_\mu$ or Aether flow $u_\mu$) acquiring a non-zero vacuum expectation value \cite{Bluhm:2004ep, Jacobson:2000xp}. The generalized $m$-Kropina metric, setting $\beta = B_\mu y^\mu$, provides a geometrical realization of the kinematics for test particles navigating this background field.  Our rigidity theorems impose severe physical constraints on these models. For $m \in 2\mathbb{Z}$. If Einstein generalized $m$-Kropina spacetime with $m$ is fractional,  by Theorem \ref{main results Einstein}, the spacetime must be strictly Ricci-flat ($Ric = 0, R_{ij} = 0$). In the context of the Finslerian field equations, this wipes out the standard purely geometric curvature terms, meaning any dynamical solution must be governed entirely by the non-Riemannian Landsberg terms ($J_i$). It severely restricts the geometry, dictating that the gravitational dynamics cannot couple to standard isotropic curvature, but must arise purely from the anisotropic interaction with the background vector field. Finally, the exact rational expressions for the geodesic spray coefficients ($G^\mu$) provide a direct kinematic pathway to constrain SME coefficients via astrophysical observations.
\subsubsection{Astrophysical and Analogue Limits: Ergospheres and Quantum Fluids}
In relativistic astrophysics, evaluating the generalized $m$-Kropina metric at $m=1$ yields $F = c \frac{\alpha^2}{\beta} + r\beta$, representing a Randers deformation of a Kropina metric. The pure Kropina term ($c\alpha^2/\beta$) famously governs critical boundary phenomena, such as the stationary limit surface of a rotating black hole's ergosphere, where the frame-dragging ``wind'', see e.g. \cite{Erasmo}, reaches the speed of light and the temporal Killing vector becomes null \cite{Gibbons:2008zi}. The linear Randers 1-form ($r\beta$) dynamically shifts this boundary, encoding an additional physical drift induced by the Lorentz-violating background. 

Beyond astrophysics, this geometric framework directly translates to analogue gravity, where acoustic perturbations in quantum fluids (such as Bose-Einstein Condensates) propagate through an effective curved spacetime \cite{Barcelo:2005fc}. Because our generalized metric intrinsically incorporates fractional anisotropy, it could generate Modified Dispersion Relations (MDRs) purely from the underlying geometry.  This suggests a possible geometric framework for studying modified dispersion relations in analogue-gravity settings.
\subsubsection{An anisotropic flat spacetime and Very Special Relativity (VSR)}
 Let $M$ be $4$-dimensional manifold with coordinates $(x^1, x^2 , x^3 ,  x^4 )$.
 An $m$-Kropina metric $ F= \frac{\alpha ^{m+1}}{\beta ^m}$,  with $m\neq 0,-1$, defined on $M$ by
    \begin{equation} \label{eq:unique_nonnull_ricci_flat_mKrop}
        F = \left| \eta_{ij}dx^i d x^j\right|^{(1+m)/2}\left(c_i d x^i\right)^{-m},
    \end{equation}
    where $\eta_{11}=-1,\, \eta_{22}=\eta_{33}=\eta_{44}=1$,  other values of $\eta_{ij}=0$ and $c_i$'s are real constants. It is clear that $\beta$ has constant length, given by,
    $$b^2 :=||\beta ||^2_{\alpha} = \eta^{ij} c_i c_j = -c_1^2 + c_2^2+c_3^2+c_4^2.$$
    This $F$ with  $b \neq 0$ is Ricci-flat and $(M,F)$ represents flat spacetime with anisotropic structure model. It can be seen from \cite[Proposition 6.4.1]{Sjors thesis},  a $4$-dimensional  $m$-Kropina space with $||\beta ||^2_{\alpha}\neq 0$ is Ricci-flat if and only if $\alpha  =  \left| \eta_{ij}dx^i d x^j\right|^{1/2}$ and $\beta = c_i d x^i $. 
    
 If $F$, defined by \eqref{eq:unique_nonnull_ricci_flat_mKrop}, with $b=0$, we get the Very Special Relativity (VSR) and the Bogoslovsky spacetime in which Lorentz symmetry is restricted to a subgroup preserving a preferred null direction. The kinematic realization of this is the Bogoslovsky spacetime \cite{Cohen:2006ky, Bogoslovsky:1977}. 
    
    Now,  it is clear that the pseudo-Riemannian metric $\alpha = \eta_{ij}dx^i d x^j$ is Ricci-flat and the $1$-form $\beta$ is parallel with respect to $\alpha$ (as the covariant derivative of $\beta$ w.r.to $\alpha$ vanishes, $b_{i ; j} =0$).
   Moreover,  the Finsler metric \eqref{eq:unique_nonnull_ricci_flat_mKrop} represents an exact  solution of each of Finslerian field equation \eqref{general field eq}.
   \begin{Rem}
The metric \eqref{eq:unique_nonnull_ricci_flat_mKrop} leads to  the following generalized $m$-Kropina metric  with the same geometric properties:
 $$F=\left(c\left| \eta_{ij}dx^i d x^j \right|+r c_i c_j d x^i dx^j \right)^{\frac{1+m}{2}} \, \left(c_i d x^i\right)^{-m},  $$ on the conic domain where   $ rs^2 \neq cm. $ 
 \end{Rem}
\subsubsection{Finsler VSI spacetimes metric}
 Let $M$ be a $4$-dimensional smooth manifold with coordinates $(u,v,x,y)$. The following metric is a Ricci-flat $m$-Kropina Finsler metric of Berwald type with $m\neq -1,0,1$ and given by \cite{Sjors m-Kropina} 
\small{
\begin{equation*}\label{Finsler VSI}
F = \left|-2d u\, d v +\Big{(} \tfrac{3m-1-m^2}{12(m-1)^2}\, x^4 + \tfrac{1}{12} y^4+ x\, v \Big{)}du^2 + 2x\,y\, d u\, d y + d x^2 + d y^2 \right|^{\frac{1+m}{2}}\left(d u\right)^{-m}.
\end{equation*}}
Clear that the $1$-form  $\beta=d u$ is null, that is, $b= 0$.  This metric name  comes from the property of its pseudo-Riemannian part $\alpha$  which has vanishing scalar curvature invariants (VSI), see for more details \cite{VSI-4D}.  Moreover,  this $F$ represents an exact solution of the Finslerian field equation \eqref{general field eq}.
\subsubsection{A Finslerian cosmological model of spacetime in Finsler gravity}
 Let $M$ be $4$-dimensional manifold with coordinates $(t, x , y ,  z )$. Define $F$ by
\begin{align}\label{eq:metricex}
    F &= |d t^2 - A^2(t)\left(d x^2+d y^2+d z^2\right)|^\frac{m+1}{2}\left(c\, A(t)^{\tfrac{1}{m}}\,d t\right)^{-m}
\end{align}
where $c$ is a positive constant and $A(t)$ is an arbitrary function.  One can show that the $m$-Kropina metric $F$ defined by \eqref{eq:metricex} is a Ricci-flat Berwald metric.  Therefore, $F$ represents a solution to the field equation in Finsler gravity \eqref{general field eq}  and it gives  a cosmological model of spacetime in the context of Finsler gravity \cite{Fuster:2018djw}.  The anisotropic field $\beta = c\, A(t)^{\tfrac{1}{m}}$ modifies the geometry. Thus the vacuum Finsler equation admit solutions with exponentially expanding ($A(t)= A_0 \, exp(Ht),\, t>0$) or power-law scale factors ($A(t)=A_0 \, t^l, \, t>0$) without introducing explicit matter or dark-energy source terms.

\begin{center}
\textbf{Conclusion}
\end{center}


\par By proving that generalized $m$-Kropina metrics belong to the class of Almost Rational (AR) Finsler metrics, we guaranteed rationality of the geodesic spray coefficients.  We provided a table which summarizes the rationality conditions for the fundamental geometric objects derived in this work under some conditions on the anisotropic degree $m$.

Furthermore, our geometric derivations revealed profound structural constraints on the spacetime manifold dictated by $m$, each carrying some physical implications: we established a rigidity theorem demonstrating that an Einstein $m$-Kropina metric with a fractional degree ($m \notin \mathbb{Z}$) must be intrinsically Ricci-flat ($Ric = 0$).  Also,  we demonstrated under $m$ even-integer constraint, assumptions of isotropic mean Berwald or Landsberg curvatures naturally force the spacetime to be weakly Berwaldian or Landsbergian, respectively.  A Berwaldian reduction is compatible with a geometric realization of the weak equivalence principle (WEP). Finally, by applying these constraints to four-dimensional manifolds, we successfully generated conditions for getting exact solutions to the Pfeifer-Wohlfarth vacuum field equations. Our case studies, ranging from the Bogoslovsky spacetime of Very Special Relativity (VSR) to Bumblebee gravity models  illustrate the potential of generalized $m$-Kropina metrics to geometrically realize several phenomenological models purely through the spacetime fabric.  Future research can seamlessly map these four-dimensional geometries to physical observables. Prime targets for observational tests include the derivation of the parameterized post-Newtonian (PPN) limits,  for example, the calculation of orbital perihelion precession in the presence of a strong background $\beta$ field, and the study of black hole shadows in strongly Lorentz-violating regimes.
\begin{center}
\textbf{Acknowledgement} 
\end{center} 
The author wishes to express her deepest gratitude to Professor Erasmo Caponio for his insightful comments, discussions and suggestions, which greatly enhanced this work, especially conversation about the Cheng-Shen Finslerian field equation for positive-definite generalized \(m\)-Kropina metrics.

\end{document}